\documentclass{amsart}
\usepackage{amssymb}
\usepackage{amsthm}
\usepackage{amscd}
\usepackage{graphics}

\newtheorem*{Theorem 1}{Theorem 1}
\newtheorem*{Theorem 2}{Theorem 2}
\newtheorem*{sweepout}{Sweepout Estimate}
\newtheorem*{isolemma}{Isoperimetric Lemma}
\newtheorem*{fedflem}{Federer-Fleming Isoperimetric Inequality}
\newtheorem{estimate}{Estimate}
\newtheorem*{Estimate 1}{Estimate 1}
\newtheorem*{Estimate 2}{Estimate 2}
\newtheorem*{key lemma}{Key Lemma}
\newtheorem{theorem}{Theorem}
\newtheorem{lemma}{Lemma}[section]
\newtheorem{corollary}{Corollary}[section]
\newtheorem*{cor}{Corollary}
\newtheorem{prop}{Proposition}[section]

\numberwithin{equation}{section}
\title{Area-expanding embeddings of rectangles}
\author{Larry Guth}
\address{Department of Mathematics, Stanford, Stanford CA, 94305 USA}
\email{lguth@math.stanford.edu}

\begin{document}

\begin{abstract} We estimate whether there is a k-expanding
embedding from one n-dimensional rectangle into another.  Our
estimates are accurate up to a constant factor $C(n)$.
\end{abstract}

\maketitle

Suppose that $U, V \subset \mathbb{R}^n$ are open sets.  An
embedding $I: V \rightarrow U$ is called k-expanding if, for
every k-dimensional surface $\Sigma \subset V$, the volume of
$I(\Sigma)$ is at least the volume of $\Sigma$.  Our theorem
describes when there is a k-expanding embedding from one
n-dimensional rectangle into another.  It is sharp up to a
constant factor in each dimension.

\begin{theorem} For each dimension $n$, there is a constant $c(n)
> 0$ so that the following holds.  Let $R$ be an n-dimensional
rectangle with dimensions $R_1 \le ... \le R_n$, and let $S$ be
an n-dimensional rectangle with dimensions $S_1
\le ... \le S_n$.

If there is a k-expanding embedding from $S$ into $R$, then,
for all integers $j, l$ in the ranges $0 \le j < k \le
l \le n$,

$$ (R_1 ... R_j)^{\frac{l-j}{k-j}} R_{j+1} ... R_l \ge c(n) (S_1
... S_j)^{\frac{l-j}{k-j}} S_{j+1} ... S_l. \eqno{(*)}$$
\end{theorem}

\begin{theorem}

Conversely, for each dimension $n$ there is a constant $C(n) > 0$
so that the following holds.  If, for all integers $j, l$ in the
ranges $0 \le j < k \le l \le n$, 

$$ (R_1 ... R_j)^{\frac{l-j}{k-j}} R_{j+1} ... R_l \ge C(n) (S_1
... S_j)^{\frac{l-j}{k-j}} S_{j+1} ... S_l, \eqno{(**)}$$

\noindent then there is a k-expanding embedding from $S$ into $R$.

\end{theorem}

Note that the necessary conditions $(*)$ and the sufficient
conditions $(**)$ are identical except that the constant $c(n)$
is replaced by the larger constant $C(n)$.

Some special cases of Theorem 1 were proven in \cite{G1}.  The
main contribution of this paper is to prove Theorem 1 in the
remaining harder cases.  In order to put the new methods in
context, we give an overview of the problem, starting with
the simplest cases.

\vskip5pt

{\bf Overview of area-expanding embeddings}

\vskip5pt

We begin by discussing the two easy cases $k=n$ and $k=1$.  If
$k=n$, then $(*)$ reduces to the one inequality $R_1 ... R_n
\gtrsim S_1 ... S_n$, which says that the volume of $R$ is bigger
than the volume of $S$.  This one condition is sufficient for
finding an n-expanding embedding from $S$ into $R$.  For example,
one can find a linear n-expanding embedding.

If $k=1$, then $(*)$ says $R_1 ... R_l \gtrsim S_1 ... S_l$ for
each $1 \le l \le n$.  These inequalities say that the smallest
l-dimensional cross-section of $S$ has smaller volume than the
smallest l-dimensional cross-section of $R$.  We'll say more
about the proof of this inequality a little lower in the
introduction.  In order to prove Theorem 2, we need to use
nonlinear maps.  For example, suppose that $R$ is the unit square
and that $S$ is a long thin rectangle with dimensions $(1/2)
\epsilon \times (1/2) \epsilon^{-1}$, $\epsilon < 1/10$.  There
is no linear 1-expanding embedding from $S$ into $R$, but there
is a non-linear 1-expanding embedding that folds $S$ into $R$, as
shown in the following figure.

\vskip5pt

\includegraphics{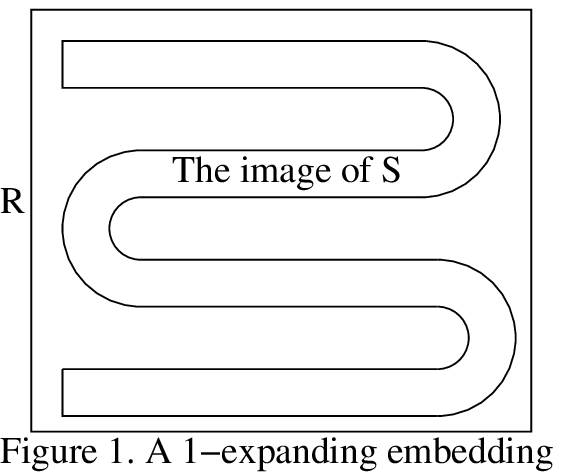}

\vskip5pt

\noindent Using these folding maps repeatedly, it's not hard to
prove Theorem 2 for $k=1$.

With that background, we turn to the main case $2 \le k \le
n-1$.  To construct k-expanding embeddings, we use the two methods
above.  We use k-expanding linear maps, and we also use simple
folding maps like the one in Figure 1.  Composing these two kinds
of maps, we construct enough embeddings to prove Theorem 2.

We now return to the proof of Theorem 1, which is the main subject
of the paper.  For general $k$, the inequalities in $(*)$ can be
divided into two kinds.  First, we have the inequalities $R_1 ...
R_l \gtrsim S_1 ... S_l$ for each $k \le l \le n$.  We have
already seen this kind of inequality about cross-sectional
volumes in the case $k=1$.  Second, we have more complicated
inequalities with $j > 0$.  These more complicated inequalities
appear only when $k$ is in the range $2 \le k \le n-1$.  
For example, if $k=2$, we have the
inequality $R_1^2 R_2 R_3 \gtrsim S_1^2 S_2 S_3$.

The proof
of the first inequalities $R_1 ... R_l \gtrsim S_1 ... S_l$
follows from a sweepout estimate as follows.  The rectangle $R$ 
may be sliced into parallel l-dimensional
rectangles with dimensions $R_1 \times ... \times R_l$.  If we
take the pullback of these surfaces in $S$, then we get a family
of surfaces sweeping out the rectangle $S$.  We refer to these
surfaces as slices of $S$.  This construction is
illustrated in the figure below.

\vskip5pt

\includegraphics{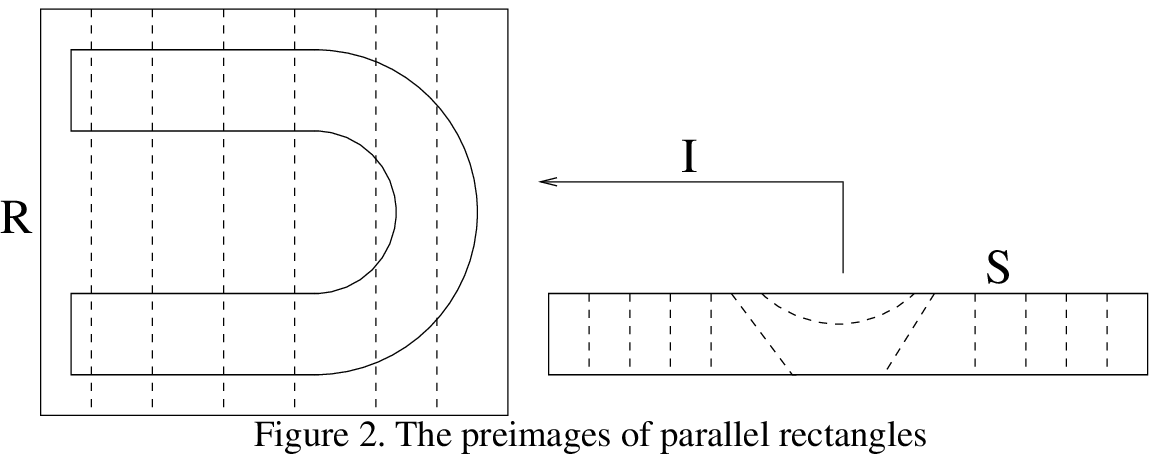}

\vskip5pt

Now $R_1 ... R_l$ is the volume of each rectangular slice of $R$. 
It follows from linear algebra that if $l \ge k$, then a
k-expanding map is also l-expanding.  (The linear algebra is
described in Appendix 1.)  Therefore, each slice of $S$ has
volume at most $R_1 ... R_l$.  Next we apply the sweepout
estimate of Almgren and Gromov.

\begin{sweepout} (Almgren, Gromov \cite{Gr}, \cite{G1}) A family of
l-dimensional surfaces sweeping out $S$ contains a surface of
volume at least $c(n) S_1 ... S_l$.
\end{sweepout}

Each slice of $S$ has volume at most $R_1 ... R_l$, but one slice
of $S$ has volume at least $c(n) S_1 ... S_l$, and so we conclude
that $R_1 ... R_l \gtrsim S_1 ... S_l$, proving $(*)$ in the case
$j=0$.

If $j > 0$, then the algebra in $(*)$ is complicated looking.  We
can think of $(*)$ as a statement about the j-dimensional
width and the l-dimensional width of $R$ and $S$.  This point of
view becomes clearer if we rewrite $(*)$ in the following
equivalent way.

$$R_1 ... R_l \gtrsim [(S_1 ... S_j) / (R_1 ...
R_j)]^{\frac{l-k}{k-j}} S_1 ... S_l.$$

\noindent If there is a k-expanding embedding from $S$ into $R$, we already
know that $R_1 ... R_l \gtrsim S_1 ... S_l$.  If $R_1 ... R_j \ge
S_1 ... S_j$, then $(*)$ follows automatically.  So we only need
to consider the case that $R_1 ... R_j$ is much smaller than $S_1
... S_j$.  In this case, $(*)$ says that the l-dimensional width
of $R$ must be substantially larger than the l-dimensional width
of $S$: larger by a factor $\sim [(S_1 ... S_j) / (R_1 ...
R_j)]^{\frac{l-k}{k-j}}$.  In other words, it is possible to
squeeze $S$ into a rectangle $R$ with much smaller j-dimensional
width only if $R$ has much larger l-dimensional width.

\vskip5pt

{\bf The tightening construction}

\vskip5pt

Now we describe the new technique in this paper.  As in Figure 2,
we look at the preimages in $S$ of parallel l-dimensional
rectangles in $R$.  We will give a proof by contradiction, so we
assume that $(*)$ is violated.  If $j=0$, we saw above that the
slices of $S$ do not have enough volume to sweep out $S$.  If $j
> 0$, then the slices of $S$ have enough volume to sweep out $S$,
but in a subtler way, we will show that they are still not big
enough to sweep out $S$.  The rough idea is that since $[0, R_1]
\times ... \times [0, R_l]$ is shaped very differently from $[0,
S_1] \times ... \times [0, S_l]$, the slices in $S$ have to
``scrunch up''.

\vskip5pt

\includegraphics{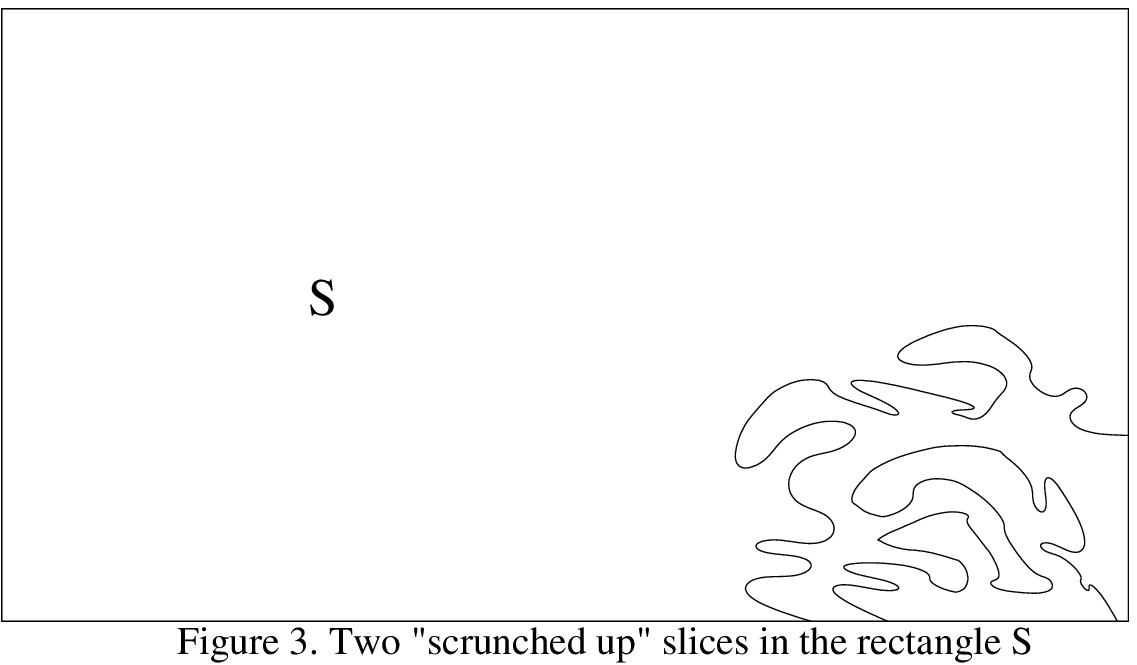}

\vskip5pt

The two curves in Figure 3 are long enough to stretch from the
bottom of $S$ to the top of $S$, but they are too scrunched up to
do so.  If a family of curves sweeps out the rectangle $S$, then
they cannot all be as scrunched up as these.

In order to prove that the slices are ``scrunched up'', and in
order to exploit this scrunching, we proceed as follows.  We
subdivide the rectangle $[0, R_1] \times ... \times [0, R_l]$
into subrectangles at a well-chosen scale.  Each slice of $S$ is
thus subdivided into pieces given by the inverse images of the
subrectangles.  A subdivision of the slices is shown in Figure 4.

\vskip5pt

\includegraphics{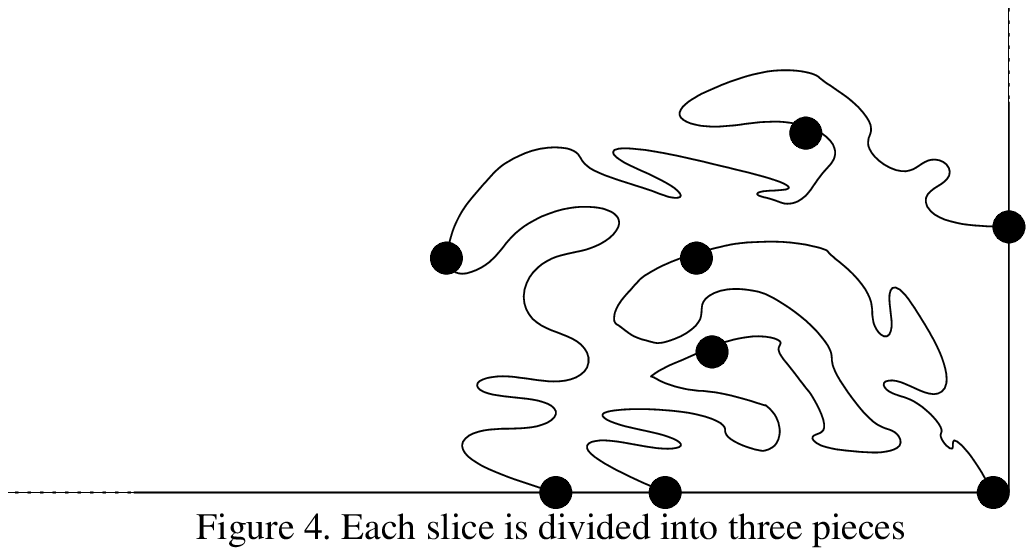}

\vskip5pt

Figure 4 shows a magnified view of the curves from Figure 3. Each
slice has been subdivided into three pieces.  The large dots mark
the endpoints of the pieces.

Now, we use an isoperimetric inequality to ``tighten'' each piece
of each slice.  This is the key step in the proof.  It involves a
new variant of the isoperimetric inequality.  We describe it in
more detail below.  Continuing informally, we show the new
tightened slices in Figure 5.

\vskip5pt

\includegraphics{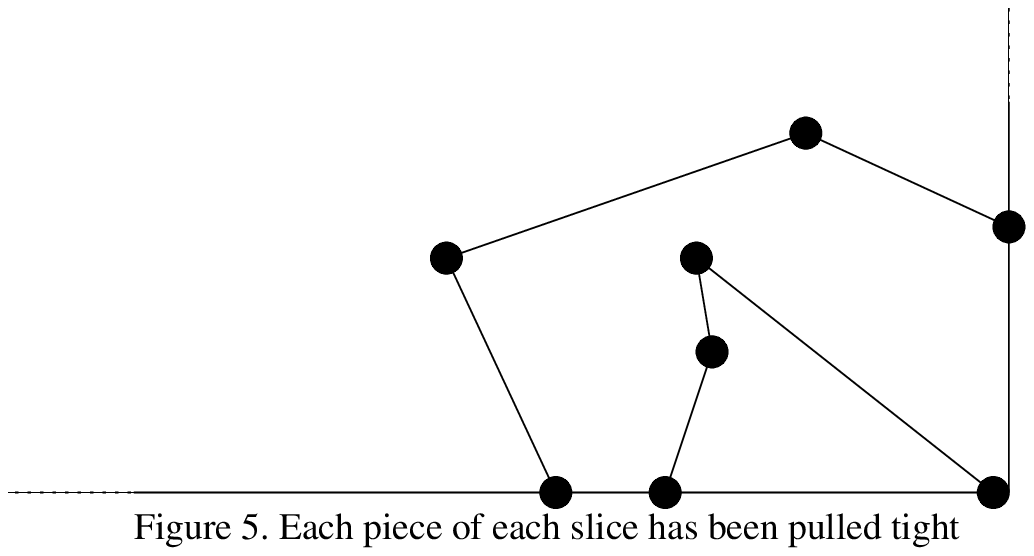}

\vskip5pt

The dots in Figure 5 are in the same locations as the dots in
Figure 4, but instead of connecting them with scrunched up curves
we have connected them with straight lines.  In the body of the
paper, the tightened pieces are not completely flat but are
constructed by a Federer-Fleming type argument.  This tightening
reduces the volume of the slices, and if $(*)$ is violated, then
the tightened slices do not have enough volume to sweep out $S$. 
This finishes our cartoon outline of the proof of Theorem 1.

The key step of tightening the pieces is done with the help of
an isoperimetric inequality.  In order to do this, we need to choose
a piece of the slice so that the piece itself has a large volume
but the boundary of the piece has a small volume.  These pieces
exist because $R_1 ... R_j$ is much smaller than $S_1 ... S_j$.

Let's give a more precise description in a simple example.  Suppose
that $j=1$, $k=2$, and $l=3$.  Furthermore, suppose that $R_1 = 1$
and that $R_2$ and $R_3$ are much bigger than 1.  Now we divide
the rectangle $[0, R_1] \times [0, R_2] \times [0, R_3]$ into
subrectangles of dimensions $1 \times L \times L$ for a large number
$L < R_2$.  One of these subrectangles has volume $L^2$.  Its relative
boundary has area $4L$.  (The relative boundary of the subrectangle
consists of four faces with dimensions $1 \times L$.  The absolute boundary
also contains two large faces with dimensions $L \times L$, but these large
faces lie in the boundary of $R$.)  For comparison, notice that a 2-cycle
$z$ in $\mathbb{R}^n$ with area $4L$ must bound a 3-chain with volume
$\lesssim L^{3/2}$, which is much smaller than $L^2$.  So the
subrectangle is a large 3-chain with a small relative boundary.

The preimage of this subrectangle in $S$ is a relative 3-chain with
a small relative boundary.  In order to tighten it, we need to prove
an isoperimetric inequality for relative cycles in the rectangle $S$.
In particular, we will prove and use the following estimate for relative
integral cycles in $S$.

\begin{isolemma} If $z$ is a p-dimensional relative cycle in
$S$ with volume $c(n) S_1 ... S_j A^{p-j}$ for some $A$ in the range
$S_j < A < S_{j+1}$, then $z$ bounds a (p+1)-chain with volume at
most $C(n) S_1 ... S_j A^{p-j+1}$.
\end{isolemma}

The Isoperimetric Lemma is a modification of the
Federer-Fleming isoperimetric inequality.  For reference, we
recall their inequality.

\begin{fedflem} If z is a q-dimensional cycle in
$\mathbb{R}^n$ with volume $A^q$, then $z$ bounds a
(q+1)-dimensional chain with volume at most $C(n) A^{q+1}$.
\end{fedflem}

The inequality in the Isoperimetric Lemma depends on the dimensions of $S$.
This is necessary: there is no isoperimetric inequality for relative
cycles that holds uniformly for all rectangles.  Instead, there is a different
isoperimetric profile for each rectangle, and we have to estimate how
the profile depends on the dimensions of the rectangle.  In the paper,
we give a fairly precise description of this isoperimetric profile,
and the Isoperimetric Lemma above is a special case.

The algebra in the Isoperimetric Lemma is somewhat complicated. 
To understand it, it helps me to consider the special case that
$z$ has the form $[0, S_1] \times ... \times [0, S_j] \times z'$,
where $z'$ is an absolute $(p-j)$-dimensional cycle in $[0,
S_{j+1}]
\times ... \times [0, S_n]$.  The cycle $z'$ would have volume
$c(n) A^{p-j}$, and the Federer-Fleming inequality implies that
$z'$ bounds a chain $y$ of volume at most $C(n) A^{p-j+1}$. 
Hence $z$ bounds $[0, S_1] \times ... \times [0, S_j] \times y$,
which has volume at most $C(n) S_1 ... S_j A^{p-j+1}$.  The
Isoperimetric Lemma says that the same estimate holds for a
general cycle $z$ as long as the volume of $z$ lies in an
appropriate range.  We prove it by using the construction of
Federer-Fleming at a sequence of different scales.

We are now ready to fill in all the details in the cartoon
outline above.  In order to keep the algebra simple, we again
focus on the special case $j=1$, $k=2$, $l=3$.  In this case,
condition $(*)$ reads $R_1^2 R_2 R_3
\gtrsim S_1^2 S_2 S_3$.  We already know that $R_1 R_2 R_3 \gtrsim
S_1 S_2 S_3$, so we only need to prove $(*)$ in the case that $R_1$ is
much smaller than $S_1$.  We consider a 3-dimensional rectangle
in $R$ with dimensions $[0, R_1] \times [0, R_2] \times [0,
R_3]$, parallel to the smallest 3-face of $R$.  We let $z$ denote
the inverse image of this rectangle in $S$.  The relative cycle
$z$ is one of the scrunched up slices in Figure 3.

Next we divide $z$ into pieces.  First we subdivide the rectangle
$[0, R_1] \times [0, R_2] \times [0, R_3]$ into subrectangles of
dimensions $R_1 \times L \times L$, for a number $L > R_1$, which
we choose later.  We let $C_i$ be the inverse images of these
subrectangles in $S$.  The chains $C_i$ are the pieces of the
slices in Figure 4.  We have $z = \sum C_i$, and we know that
each chain $C_i$ has volume at most $R_1 L^2$. 

Now we look at the boundaries of the chains $C_i$.  Each of our
3-dimensional subrectangles of dimension $R_1 \times L \times L$
has a relative boundary with area at most $4
R_1 L$.  Since the map $I$ is 2-expanding, the relative boundary
of each chain $C_i$ has area at most $4 R_1 L$.

Now we apply the Isoperimetric Lemma to the boundary of $C_i$. 
To make the proof work, we have to choose $L$ so that $R_1 L$ is
between $S_1^2$ and $S_1 S_2$.  Then the Isoperimetric Lemma
guarantees that $\partial C_i$ bounds some 3-chain $C_i'$ with
volume at most $\sim (R_1 / S_1) R_1 L^2$.  In other words, our
bound for the volume of $C_i'$ is better than the bound for the
volume of $C_i$ by a factor $\sim (R_1 / S_1)$.  To tighten the
slice $z$, we replace each chain $C_i$ with the chain $C_i'$. 
The chains $C_i'$ are the segments in Figure 5.  We define a
relative cycle $z' = \sum C_i'$.  The relative cycle $z'$ is one
of the tightened slices in Figure 5.  The 
total volume of $z'$ is at most
$\sim (R_1 / S_1) R_1 R_2 R_3$. 

We perform the same tightening operation on every slice.  Each
tightened slice has volume at most $\sim (R_1 / S_1) R_1 R_2
R_3$.  Because of the sweepout lemma, one of the tightened slices
must have volume at least $c(n) S_1 S_2 S_3$.  Hence $(R_1 /
S_1) R_1 R_2 R_3 \gtrsim S_1 S_2 S_3$, and rearranging we get
$R_1^2 R_2 R_3 \gtrsim S_1^2 S_2 S_3$ as desired.

In general the tightening procedure is a little bit more involved.  We use
our control of the k-skeleton of the slice to tighten the
(k+1)-skeleton.  Then we use our improved control of the
(k+1)-skeleton to tighten the (k+2)-skeleton, and so on until we
get to the l-skeleton of the slice.

\vskip5pt

{\bf Complexes of cycles}

\vskip5pt

Lastly, I want to say a word about the language we use in the
proof.  We outlined our argument informally in terms of families
of surfaces, but families of surfaces are not a convenient language. 
For one problem, the tightening construction we just described
does not depend continuously on the surface.  Instead, we use a
discrete analogue of a family of cycles that we call a complex of
cycles.  Over several papers, I have found complexes of cycles to
be a simple, convenient language for arguments about
area-contracting maps.

A complex of cycles $C$ in $S$ is parametrized by a polyhedral
complex $X$.  For each p-face $F$ of $S$, the complex $C$ associates
a p-dimensional relative chain $C(F)$, and these chains are
required to fit together in a coherent way.  Figure 6 shows an
example of a complex of cycles, illustrating the way the chains
should fit together.

\vskip5pt

\includegraphics{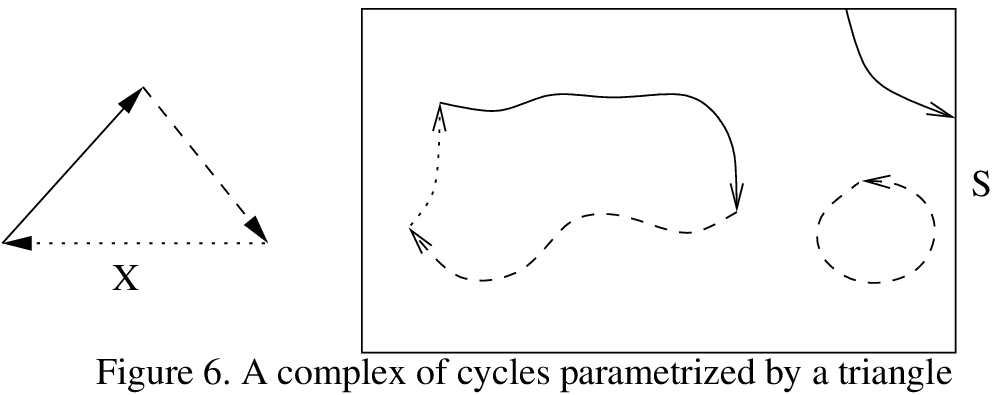}

\vskip5pt

In this example, the polyhedral complex $X$ is a triangle.  Each
side of the triangle corresponds to an oriented relative 1-chain
in $S$.  The solid line in the triangle corresonds to the two
solid curves in $S$, and so on.

Complexes of cycles were introduced by Almgren in his thesis
\cite{A} on the homotopy groups of spaces of cycles.  He begins
with a continuous family of cycles, but the first step in his
argument is to replace the continuous family by a complex of
cycles that approximates it.  Complexes of cycles were then used
by Gromov in his proof of the Sweepout Estimate \cite{Gr}.  The
first step in Gromov's proof is also to replace the continuous
family by a complex of cycles approximating it.  Almgren and
Gromov did not name the object that they use.  The name complex
of cycles comes from \cite{G4}.

\vskip5pt

Here is an outline of the paper.  In Section 1, we prove estimates 
for the isoperimetric profile of
a rectangle.  In Section 2, we state a generalization of Theorem
1.  In Section 3, we define complexes of cycles.  In Section 4 we
prove a version of the sweepout lemma for complexes of cycles. 
With this lemma, we prove Theorem 1 in the easy case $j=0$. In
Section 5, we give some algebraic preliminaries which reduce the
general case of Theorem 1 to a slightly more special case.  In
Section 6, we explain the tightening construction and prove
Theorem 1. This section is the heart of the paper.  In Section 7,
we construct area-expanding embeddings of rectangles, proving
Theorem 2.  The paper ends with two appendices.  The first
appendix covers the linear algebra related to area-expanding or
area-contracting maps.  The second appendix covers
generalizations of our results to shapes other than rectangles.

\vskip5pt

{\bf Acknowledgements.} This paper is a simplified version of the main
result of my thesis \cite{G3}.  The proof in my thesis was very
convoluted.  I am grateful to my thesis advisor, Tom Mrowka, for
his support and encouragement.

\section{The isoperimetric profile of a rectangle}

Let $R$ denote the n-dimensional rectangle $[0, R_1] \times ...
\times [0, R_n]$, where the dimensions are ordered so that $R_1
\le ... \le R_n$. In this section, we estimate the isoperimetric
profile for relative integral cycles in $R$.  Our goal is to
understand the way that the isoperimetric profile depends on the
dimensions $R_i$.

If $z$ is a relative integral k-cycle in $R$, the filling volume
of $z$ is the smallest volume of any relative (k+1)-chain $y$
with $\partial y = z$.  Let $I^k_R(V)$ denote the largest filling
volume of any k-dimensional relative integral cycle in $R$ with
volume at most $V$.

Remark: We use the following definition for volume.  If a chain
$C$ is given by $\sum c_i f_i$ where $c_i \in \mathbb{Z}$ and
$f_i$ is a Lipschitz map from the standard k-simplex to $S$, then
the volume of the chain $C$ is defined to be $\sum |c_i|
Vol(f_i^* Euc)$, where $Euc$ denotes the Euclidean metric on $R$. 
This quantity is also called the mass of $C$.  We denote the
volume of $C$ by $|C|$.

The following theorem estimates the isoperimetric profile $I^k_R$ for
the rectangle $R$.

\begin{theorem} There are constants $c(n) > 0, C(n)$ so that the
following holds.

If $V \le c(n) R_1 ... R_k$, then write $V = c(n) R_1 ... R_j
\rho^{k-j}$ for some $0 \le j \le k-1$ and some $\rho$ in the
range $R_j \le \rho \le R_{j+1}$. (These conditions determine $j$
and $\rho$ uniquely.)

$$\textrm{Then } I^k_R(V) \le C(n) R_1 ... R_j \rho^{k-j+1}.$$

In any case, $I^k_R(V) \le C(n) R_{k+1} V$.

\end{theorem}

Before we prove the theorem, we consider two examples of relative
cycles in $R$.  These examples show that our upper bounds for
$I^k_R$ are fairly sharp.  They also help me to remember the
formulas.

Pick an integer $j$ in the range $0 \le j \le k-1$.  Then
consider the cycle $[0, R_1] \times ... \times [0, R_j]
\times S^{k-j}(\rho)$ for $R_j \le \rho
\le (1/10) R_{j+1}$.  In this equation, $S^{k-j}(\rho)$ denotes a
sphere of dimension $k-j$ and radius $\rho$ contained in $[0,
R_{j+1}] \times ... \times [0, R_n]$, with center at the center
of the rectangle, $(R_{j+1}/2, ..., R_n)$.  This cycle has volume
$V \sim R_1 ... R_j \rho^{k-j}$.  The best filling of the cycle
is just $[0, R_1] \times ... \times [0, R_j]
\times B^{k-j+1}(\rho)$.  To clarify the notation, $B^{k-j+1}$ is
a Euclidean ball of dimension $k-j+1$ with boundary
$S^{k-j}(\rho)$.  This filling has volume $\sim R_1... R_j
\rho^{k-j+1}$.

Second, consider the relative cycle $[0,R_1]
\times ... \times [0, R_k] \times \{ p \}$ with multiplicity $M$, 
where $p$ is the center of the rectangle $[0, R_{k+1}] \times ...
\times [0, R_n]$.  (Alternatively, consider $M$ nearby parallel
rectangles.)  The volume of this cycle is $V = M R_1 ... R_k$. 
This cycle has filling volume $\sim M R_1 ... R_{k+1} = R_{k+1}
V$.

Remarks: These examples give lower bounds for $I^k_R(V)$.  The
lower bounds match the upper bounds in the theorem up to a
constant factor except in the delicate range $c(n) R_1 ... R_k
\le V \le R_1 ... R_k$. If $R_{k+1} >> R_k$, then the function
$I^k_R(V)$ grows very rapidly over the course of this range.  It
appears plausible that $I^k_R(V)$ is discontinuous, perhaps at the
value $V = R_1 ... R_k$.

Now we turn to the proof of the theorem.

\proof We begin by using the deformation theorem of Federer and
Fleming, which we record as a lemma.

\begin{lemma} Suppose that $z$ is a relative k-cycle in $R$.
Consider a rectangular lattice inside $R$ with each side-length
roughly equal to $L$ (up to a factor of 2), for some $L \le R_1$,
and suppose that the boundary of $R$ lies in the (n-1)-skeleton
of the lattice.  Then there is another relative cycle $z'$ in $R$
contained in the k-skeleton of the lattice and obeying the 
following inequalities.

1. The volume of $z'$ is at most $C(n) |z|$.

2. The filling volume of $z' - z$ is at most $C(n) L |z|$.

\end{lemma}

Remark: Morally, we are using a cubical lattice.  We allow a slightly
non-cubical lattice so that we can arrange for the boundary of $R$ 
to lie in the (n-1)-skeleton of the lattice.  

\proof (sketch) We sketch the proof of Federer and Fleming.  For
more details, see \cite{G3}.  We begin with a relative cycle $z$
with boundary $\partial z$ contained in $\partial R$.  We build a
sequence of homologies $z = z_n \sim z_{n-1} \sim ... \sim z_k =
z'$, where each $z_p$ has the same boundary as $z$ and $z_p$ lies
in the union of the p-skeleton of our lattice and $\partial R$. 
(If we think of $z_p$ as relative chains, then they are all
cycles and $z_p$ lies in the p-skeleton of our lattice.)

The homology from $z_p$ to $z_{p-1}$ is constructed as follows. 
For each interior p-face of our lattice, we push $z_p \cap F$
into $\partial F$ while keeping $z_p \cap \partial F$ fixed.  To
do this, we pick a random point $x$ in $F$ and push $F - \{ x \}$
radially into the boundary of $F$. For a random point $x$, this
operation stretches volume by at most a constant $C(n)$. 
Therefore, the volume of $z_{p-1}$ is at most $C(n) |z_p|$. 
Similarly, the volume of the homology from $z_p$ to $z_{p-1}$ is
at most $C(n) L |z_p|$. \endproof

Using this lemma, we prove the isoperimetric inequality by
induction on $k$. When $k=0$, $z$ is just a weighted sum of
points $\sum c_i p(i)$, where $c_i \in \mathbb{Z}$ and $p(i)$ is
a point in $R$.  The volume of $z$ is defined to be $\sum |c_i|$. 
A point $p$ with coordinates $(p_1, ..., p_n)$ bounds a segment
$[0, p_1] \times \{ p_2 \} \times ... \times \{ p_n \}$. 
Applying this operation to each point $p(i)$ with multiplicity
$c_i$, we get a filling of $z$ with volume at most $R_1 Vol(z)$.
This argument gives the base for our induction.

Now we come to the inductive step.  Suppose that $z$ is a k-cycle
with volume $V$.  We proceed in two cases.  If $V \le c(n)
R_1^k$, then we select $L = C(n) V^{1/k} \le R_1$ and pick a
rectangular lattice with sidelengths roughly $L$ and with
$\partial R$ in the (n-1)-skeleton of the lattice. Then we use
Lemma 1.1 to move $z$ to a new relative cycle $z'$ with volume at
most $C(n) V$ lying in the k-skeleton of our lattice. Since $C(n)
V \le L^k$, the new cycle $z'$ is simply $0$.  Lemma 1.1 also
guarantees us a homology from $z$ to $z'$ with volume at most
$C(n) L V$, which is at most $C(n) V^{\frac{k+1}{k}}$.  This
upper bound is the one we needed to prove.

In the second case, we suppose that $V \ge c(n) R_1^k$.  In this
case, we select $L = R_1$ and pick a rectangular lattice with
sidelengths roughly $L$ and with $\partial R$ in the
(n-1)-skeleton of the lattice.  We pick the lattice so that each
lattice point has $x_1$ coordinate either $0$ or $R_1$. Then we
use Lemma 1.1 to move $z$ to a new relative cycle $z'$ with
volume at most $C(n) V$ lying in the k-skeleton of our lattice. 
The homology from $z$ to $z'$ has volume at most $C(n) R_1 V$. 
The cycle $z'$ need not be $0$, but it is a union of interior
k-faces of our lattice.  Each interior k-face has the form $[0,
R_1] \times ...$, and so the cycle $z'$ has the special form $z'
= [0, R_1] \times z_1$ for some relative cycle $z_1$ in the
(n-1)-dimensional rectangle $[0, R_2] \times ... \times [0,
R_n]$.  The cycle $z_1$ has volume at most $C(n) V / R_1$.

By induction, we can assume that our theorem holds for $z_1$. 
Therefore, $z_1$ bounds a relative chain $C_1$ with a certain
volume bound that we calculate below.  Then $z'$ bounds $[0, R_1]
\times C_1$.  We will calculate that the volume of this filling
obeys the inequality stated in the theorem.

If the volume of $z$ is at most $c(n) R_1 ... R_k$, then the
volume of $z_1$ is at most $c(n-1) R_2 ... R_k$.  If the volume
of $z$ is equal to $c(n) R_1 ... R_j \rho^{k-j}$ for some $\rho$ in
the range $R_j \le \rho \le R_{j+1}$, then the volume of $z_1$ is
roughly $R_2 ... R_j \rho^{(k-1) - (j-1)}$ for the same $\rho$. 
By induction, $z_1$ bounds a chain $C_1$ with volume at most
$C(n-1) R_2 ... R_j \rho^{k-j+1}$, and so $[0, R_1] \times C_1$
has volume at most $C(n-1) R_1 ... R_j \rho^{k-j+1}$.  Also, the
homology from $z$ to $z'$ has volume at most $C(n) R_1 (R_1 ...
R_j) \rho^{(k-j)} \le C(n) R_1 ... R_j \rho^{k-j+1}$. Therefore,
the filling volume of $z$ is at most $C(n) R_1 ... R_j
\rho^{k-j+1}$.

In any case, $z_1$ bounds a k-chain $C_1$ of volume at most
$C(n-1) R_{k+1} V / R_1$. Hence $z' = [0,R_1] \times z_1$ bounds
a (k+1)-chain of volume at most $C(n-1) R_{k+1} V$.  The homology
from $z$ to $z'$ has volume at most $C(n) R_1 V$.  Therefore, the
filling volume of $z$ is at most $C(n) R_{k+1} V$. \endproof

The algebra above is a bit complicated.  In the sequel, we only
use the following special case, which is easier to remember.

If $z$ is a relative p-cycle in $R$ with volume $V$ at most $c(n)
R_1 ... R_j R_j^{p-j}$, then it bounds a relative (p+1)-chain $y$
in $R$ with volume at most $C(n) R_j V$.

\section{Statement of the main inequalities}

In the paper, we will prove an estimate which is a little more
general than Theorem 1.  We now formulate it in terms of
k-dilation.  Recall that the k-dilation of a smooth map $\Phi$ is
defined to be $\| \Lambda^k d \Phi \|_{L^\infty}$.  The
k-dilation measures by what factor the map $\Phi$ stretches
k-dimensional areas.  The k-dilation of $\Phi$ is at most
$\Lambda$ if and only if $\Phi$ maps every k-dimensional surface
of volume $V$ to an image of volume at most $\Lambda V$.

Recall that $R$ is an n-dimensional rectangle with dimensions
$R_1 \le ... \le R_n$ and $S$ is an n-dimensional rectangle with
dimensions $S_1 \le ... \le S_n$.  We let $Q_i$ denote the
quotient $S_i / R_i$.  We now state the main estimates of the
paper.

\begin{estimate} Suppose that $U$ is an open set in $R$ and that
$\Phi$ is a map of pairs $(U, \partial U) \rightarrow (S,
\partial S)$ of degree $D > 0$.  Suppose that $j$ and $l$ lie in
the ranges $0 \le j < k \le l \le n$.  Then the k-dilation of
$\Phi$ is bounded below by the following inequality.

$$dil_k(\Phi) \ge c(n) Q_1 ... Q_j (Q_{j+1} ...
Q_l)^{\frac{k-j}{l-j}}.$$

\end{estimate}

For example, if $I$ is a k-expanding embedding from $S$ into $R$,
then we take $U$ to be the image of $S$, and we take $\Phi$ to be
the inverse of $I$.  The map $\Phi$ has k-dilation at most 1, and
it has degree 1, and so Estimate 1 implies Theorem 1.  Estimate 1
is slightly more general because $\Phi$ need not be a
diffeomorphism. 

If the degree $D$ is large, then we can strengthen some of the
lower bounds in Estimate 1 as follows.

\begin{estimate} With the same assumptions as above, for any $0
\le j < k$, the k-dilation of $\Phi$ is bounded below by the
following inequality.

$$dil_k(\Phi) \ge c(n) D^{\frac{k-j}{n-j}} Q_1 ... Q_j (Q_{j+1}
... Q_n)^{\frac{k-j}{n-j}}.$$

\end{estimate}

In the paper \cite{G1}, I proved Estimate 1 if either $j=0$ or
$l=n$.  We will prove all the cases of Estimate 1 in this paper. 
The proof of the case $j=0$ is essentially the same as the one in
\cite{G1}, but this paper gives a new proof for the case $l=n$.

\section{Complexes of cycles}

We introduce some vocabulary that we will use in our proof.

A complex of cycles in a rectangle $S$ is a collection of chains
of different dimensions that fit together in a coherent way.  It
consists of the following data.  There is a polyhedron $X$ which
is like a parameter space for the complex. Then there is a map
$C$ which assigns to each d-dimensional face $F^d$ of $X$ a
d-dimensional relative chain in $S$.  These chains have to fit
together so that if the boundary of $F^d$ is equal to
$\sum_{i=1}^N F^{d-1}_i$, then the boundary of the chain $C(F)$
should be $\sum_{i=1}^N C(F_i)$.  In this paper, we work with
complexes of cycles over $\mathbb{Z}$, and so all the faces and
chains in the above discussion are oriented.

More formally, the map $C$ is a chain map between two complexes. 
The first complex is generated by the faces of $X$ with integral
multiplicities and the natural boundary operations.  The second
complex is the complex of integral relative Lipschitz cycles in
$S$, which we denote $I_{rel}(S)$. 

This definition is due to Almgren.  Almgren introduced it in his
paper on the topology of the space of cycles \cite{A}.  
For more explanation of the definition, see Section 1 of \cite{G4}.

We remark that the complex $X$ may have dimension bigger than $n$.
Even if $d > n$, the definition of Lipschitz d-chain in $S$ makes
sense.

We give an example of a complex of cycles.  
If $U \subset R$ is an open set and $\Phi$ is a map
from $(U, \partial U)$ to $(S, \partial S)$, then we can define a
complex of cycles by noticing where $\Phi$ maps various chains. 
Let us fix a polyhedral structure $P$ on $R$.  For each face $F$
of this structure, we define $C_\Phi(F)$ to be $\Phi(F \cap U)$. 
The complex $C_\Phi$ sends each face $F$ contained in the
boundary of $R$ to zero, and so we can say that $C_\Phi$ is
parametrized by $(R, \partial R)$.

Since $C$ is a chain map, it induces a map on homology from
$H_*(X, \mathbb{Z})$ to $H_*(S, \partial S, \mathbb{Z})$.
In particular, if $C$ is a complex of cycles
parametrized by $(R, \partial R)$, then it induces a map from
$H_*(R, \partial R, \mathbb{Z})$ to $H_*(S, \partial S,
\mathbb{Z})$.  We define the degree of $C$ to be the degree of
this map on $H_n$.  The degree of $C_{\Phi}$ is the same as the
degree of $\Phi$.

A homotopy of complexes of cycles is a complex $C$ parametrized
by $X \times [0,1]$.  If the restriction of $C$ to $X \times \{ 0
\}$ is a complex $C_0$ and the restriction of $C$ to $X \times \{
1 \}$ is $C_1$, then we say that $C$ is a homotopy from $C_0$ to
$C_1$.  If $C_0$ and $C_1$ are homotopic, then the induced maps
on homology $H_*(X, \mathbb{Z}) \rightarrow H_*(S, \partial S,
\mathbb{Z})$ are the same.

\section{The sweepout lemma}

We now prove a lemma that says that if all the chains in a
complex are small enough then the complex is null-homotopic. 
The lemma and proof are based on an argument of Gromov from
page 134 of \cite{Gr}.

\begin{lemma} There is a constant $c(n) > 0$ so that the
following estimate holds.

Suppose that $C_0$ is a complex of cycles in $S$ parametrized by $X$. 
Suppose that for each vertex $v$ of $X$, $C_0(v)$ is equal to
$0$.  Suppose that for each p-face $F^p$ in $X$, $C_0(F^p)$ has
volume at most $c(n) S_1 ... S_p$.  Then $C_0$ is null-homotopic.
\end{lemma}

Lemma 4.1 is closely related to the Sweepout Estimate stated in
the introduction.  Gromov used this argument to prove the
sweepout estimate on page 134 of \cite{Gr}.

\proof We let $C_1$ denote the zero map.  We have to prove that
$C_0$ is homotopic to $C_1$ by constructing a homotopy $C$ between
them.  The homotopy $C$ needs to be defined on $X \times [0,1]$,
and it is already defined on $X \times \{ 0 \}$ and on $X \times
\{ 1 \}$.  We define $C$ one skeleton at a time.

We will prove inductively that we can extend $C$ to the
p-skeleton of $X \times [0,1]$ while preserving the inequality
$|C(F^p)| \le c(n) S_1 ... S_p$ for all $p \le n$.  
To start the induction, we
define $C$ on the 1-skeleton by setting $C(v \times [0,1])$ equal
to zero for each vertex $v$ of $X$.  Since $C_0(v) = 0 = C_1(v)$,
this choice is allowed and it clearly obeys our volume estimate.
By induction, we may assume that we have done the extension to
the (p-1)-skeleton of $X
\times [0,1]$.  When we extend to the p-skeleton, we have to
define $C(F^p)$ for each p-face so that $\partial C(F^p) =
C(\partial F^p)$.  By induction, $C(\partial F^p)$ is a
(p-1)-cycle in $S$ with volume at most $c(n) S_1 ... S_{p-1}$. 
According to Theorem 3, we can fill this cycle with volume at
most $c(n) S_1 ... S_{p-1} S_{p-1} \le c(n) S_1 ... S_p$.

Next we have to extend $C$ to the (n+1)-skeleton of $X$.  We
have already defined $C$ on the n-skeleton.  In particular,
$C(\partial F^{n+1})$ is a relative n-cycle in $S$ with
volume at most $c(n) S_1 ... S_n < S_1 ... S_n$.  Therefore
this n-cycle is exact.  We define $C(F^{n+1})$ to be any
(n+1)-chain with the given boundary.  Finally we extend
to the higher skeleta.  There is no obstruction to finding an
extension to the higher skeleta because $H_p(S, \partial S) = 0$
for $p \le n+1$.
\endproof

(The same proof works for a complex of cycles parametrized by
$(R, \partial R)$. In this case, we get a homotopy parametrized
by $(R \times [0,1], \partial R \times [0,1])$.)

Using this lemma, we can prove the easiest cases of Estimates 1
and 2.  These cases were first proven in \cite{G1}, but we
include them here for completeness.

\begin{prop} If $U$ is an open set in $R$ and if $\Phi: (U,
\partial U) \rightarrow (S, \partial S)$ is a map of degree $D
\not= 0$, then the k-dilation of $\Phi$ is at least $c(n) Q_1 ...
Q_k$.
\end{prop}

\proof By scaling, we may assume that $\Phi$ is k-contracting and
it then suffices to prove that $R_1 ... R_k \ge c(n) S_1 ...
S_k$.  We assume that $R_1 ... R_k < c(n) S_1 ... S_k$ and
proceed to a contradiction.

We cut $R$ into rectangular blocks which are each congruent to
$[0, R_1] \times ... \times [0, R_k] \times [0, \epsilon]^{n-k}$
for some small number $\epsilon > 0$.  All the rectangular blocks
are parallel, and they form a grid of dimension $1 \times ...
\times 1 \times (R_{k+1}/\epsilon) \times ... \times (R_n /
\epsilon)$. Now we look at the complex $C_{\Phi}$ corresponding
to this decomposition.

If $p < k$, then each p-face of our decomposition lies on the
boundary of $R$ and so is mapped to $0$.  Each k-face of our
decomposition has volume at most $R_1 ... R_k$.  For each k-face
$F^k$, $C_{\Phi}(F^k)$ has volume less than $c(n) S_1 ... S_k$,
since $\Phi$ is k-contracting.  Similarly, for $p > k$, each
p-face $F^p$ has volume at most $R_1 ... R_k
\epsilon^{p-k}$.  In Appendix 1, we prove that if $\Phi$ is
k-contracting then it is also l-contracting for each $l \ge k$.
So $C_{\Phi}(F^p)$ has volume less than $R_1 ... R_k
\epsilon^{p-k}$. If we choose $\epsilon$ small enough, then Lemma
4.1 implies that $C_{\Phi}$ is null-homotopic. In particular
$C_{\Phi}$ has degree zero.  But we have already seen that
$C_{\Phi}$ has degree $D$ which we assumed non-zero. 
\endproof

Since the map $\Phi$ is also l-contracting for all $l \ge k$, we
get the following more general proposition.

\begin{prop} If $l \ge k$, if $U$ is an open set in $R$, and if
$\Phi: (U, \partial U) \rightarrow (S, \partial S)$ is a map of
degree $D \not= 0$, then the k-dilation of $\Phi$ is at least
$c(n) (Q_1 ... Q_l)^{k/l}$.  Also, the k-dilation of $\Phi$ is at
least $(|D| Q_1 ... Q_n)^{k/n}$.
\end{prop}

\proof By the last proposition, the l-dilation of $\Phi$ is at
least $c(n) Q_1 ... Q_l$.  Also, the n-dilation of any degree $D$
map is at least $|D| Q_1 ... Q_n$.  Therefore, the k-dilation of
$\Phi$ is at least $c(n) (Q_1 ... Q_l)^{k/l}$ and at least $(|D|
Q_1 ... Q_n)^{k/n}$. \endproof

Proposition 4.2 proves Estimates 1 and 2 in the case $j=0$.

\section{Algebraic preliminaries}

We rewrite the remaining cases of our estimates.

\begin{Estimate 1} (Non-trivial cases) There is a constant $c(n) >
0$ so that the following holds.  Let $R,S$ be n-dimensional
rectangles.  Suppose $U \subset R$ is an open set.  Suppose that
$\Phi$ is a k-contracting map from $U$ to $S$ of degree $D \not=
0$.  Suppose $0 < j < k < l$.

$$\textrm{Then } [(R_1 ... R_j)/(S_1 ... S_j)]^{\frac{l-k}{k-j}}
R_1 ... R_l \ge c(n) S_1 ... S_l. \eqno{(1)}$$

\end{Estimate 1}

\begin{Estimate 2} (Non-trivial cases) In the same situation as
above, the following inequality holds.

$$[(R_1 ... R_j)/ (S_1 ... S_j)]^{\frac{n-k}{k-j}} R_1 ... R_n
\ge c(n) |D| S_1 ... S_n. \eqno{(2)}$$

\end{Estimate 2}

Fix $j$.  We define $L$ by the equation $R_1 ... R_j L^{k-j} =
\delta(n) S_1 ... S_j S_j^{k-j}$, where $\delta(n) > 0$ is a
small dimensional constant.

In the next section, we will prove Estimates 1 and 2
under the assumption that $L \le R_{j+1}$.

We now check that it suffices to prove the estimates in this
special case.  This checking just takes a little algebra - the
geometric part of our proof is in the next section.

We can rewrite our inequalities as follows.

$$[(R_1 ... R_j)/(S_1 ... S_j)]^{\frac{1}{k-j}} \ge 
c(n) [(S_1 ... S_l)/(R_1 ... R_l)]^{\frac{1}{l-k}}. \eqno{(1')}$$

$$[(R_1 ... R_j)/(S_1 ... S_j)]^{\frac{1}{k-j}} \ge 
c(n) [|D| (S_1 ... S_n)/(R_1 ... R_n)]^{\frac{1}{n-k}}. \eqno{(2')}$$

The right-hand sides of both equations is independent of $j$.  So
it suffices to pick the one value of $j$ that minimizes the
left-hand side and to prove our theorem for this one value of
$j$.  Now for this value of $j$, we will prove that $L \le
R_{j+1}$. 

We see this inequality in two cases.  In the first case, it may
happen that $j = k-1$.  In this case, $L = \delta(n) S_1 ...
S_{k-1} S_{k-1} / (R_1... R_{k-1}) \le \delta(n) S_1 ... S_k /
(R_1 ... R_{k-1})$.  But by Proposition 4.1, $R_1 ... R_k \ge
c(n) S_1 ... S_k$.  Therefore, $\delta(n) S_1 ... S_k / (R_1 ...
R_{k-1}) \le \delta(n) c(n)^{-1} R_k$. If we choose $\delta(n)$
small enough, then $\delta(n) c(n)^{-1} R_k
\le R_k = R_{j+1}$.

In the second case $j < k-1$.  In this case, $j+1$ was a legal
competitor for $j$, and so we conclude that

$$[(R_1 ... R_j)/(S_1 ... S_j)]^{\frac{1}{k-j}} \le [(R_1 ...
R_{j+1})/(S_1 ... S_{j+1})]^{\frac{1}{k-j-1}}.$$

We raise each side of the equation to the power $(k-j)(k-j-1)$
and then move all the powers of $R$ to the righthand side.

$$S_1 ... S_j S_{j+1}^{k-j} \le R_1 ... R_j R_{j+1}^{k-j}.$$

A fortiori, $S_1... S_j S_j^{k-j} \le R_1 ... R_j R_{j+1}^{k-j}$.   
On the other hand, $S_1 ... S_j S_j^{k-j} \ge R_1 ... R_j L^{k-j}$.  
Therefore, $L \le R_{j+1}$.

\section{Tightening a complex of cycles}

In this section we prove our main estimates by cutting the
rectangle $R$ into pieces, mapping the pieces into $S$, and then
pulling them tight with the isoperimetric inequality.  To begin,
we cut $R$ into subrectangles of a carefully chosen size.

We define $L$ by the equation $R_1 ... R_j L^{k-j} =
\delta(n) S_1 ... S_j S_j^{k-j}$, where $\delta(n) > 0$ is a
small constant that we can choose later.  In this section we will
make the mild assumption that $L \le R_{j+1}$.  In Section 5, we
explained how the general case follows from this special case by
high-school algebra.  We pick a polyhedral structure on $R$ by
cutting it into rectangular blocks of dimensions $R_1 \times ...
\times R_j \times L \times ... \times L$.  (By making a mild
change in the dimensions of $R$, we may also assume that $L$
divides $R_i$ for each $i \ge j+1$.)

We let $B$ be the chain complex generated by the interior faces
of this decomposition. The homology of $B$ is $H_*(R, \partial R,
\mathbb{Z})$. We let $C_0$ be the chain map $B \rightarrow
I_{rel}(S)$ associated to $\Phi$.  In other words, if $F$ is a
face of $B$, then $C_0(F)$ is $\Phi(F \cap U)$.  The degree of
$C_0$ is $D$, the degree of $\Phi$.

By repeatedly using the isoperimetric inequality, we will
``tighten'' $C_0$ to a new complex of cycles $C_1$.

The complex $C_1$ agrees with $C_0$ for faces of dimension at
most $k$.  For faces of higher dimension, $C_1$ is different from
$C_0$.  We define $C_1$ by induction on the dimension.

First we define $C_1(F^{k+1})$.  We have already defined
$C_1(\partial F^{k+1})$.  Each face of $\partial F^{k+1}$ has
k-volume at most $\delta(n) S_1 ... S_j S_j^{k-j}$, and so
$C_1(\partial F^{k+1})$ has volume at
most $\delta(n) 2^n S_1 ... S_j S_j^{k-j}$.  If we pick $\delta$
small enough, the isoperimetric inequality tells us that
$C_1(\partial F^{k+1})$ bounds a (k+1)-chain with volume at most
$\delta(n) C S_1 ... S_j S_j^{k-j+1}$.  We define $C_1(F^{k+1})$
to be a (k+1)-chain with this volume bound.  We repeat this
construction for every (k+1)-face in our decomposition of $R$.

Then we proceed inductively, defining $C_1$ one skeleton at a
time so that at each stage it obeys the inequality $|C_1(F^p)| <
\delta C S_1 ... S_j S_j^{p-j}$.  Suppose we have defined $C_1$
on the p-skeleton and that $F^{p+1}$ is a (p+1)-face.  We have
already defined $C_1(\partial F^{p+1})$ and it has volume at most
$\delta C S_1 ... S_j S_j^{p-j}$.  Assuming $\delta$ is
sufficiently small, we can apply the isoperimetric inequality
to fill $C_1(\partial F^{p+1})$ by a (p+1)-chain of
volume at most $\delta C S_1 ... S_j S_j^{p+1-j}$. We define
$C_1(F^{p+1})$ to be this chain.

A key point in the proof is that the new complex $C_1$ has the same degree
as the original complex $C_0$.

\begin{key lemma} The degree of $C_1$ is equal to the degree of $\Phi$.
\end{key lemma}

\noindent This point is the subtlest part of our argument, and so we defer
the proof until the end of the section.

\vskip5pt

{\bf Gluing a complex of cycles}

\vskip5pt

In place of $B$, we now consider a coarser decomposition of the
rectangle $R$.  This time we divide $R$ into blocks with
dimensions $R_1 \times ... \times R_l \times L \times ... \times
L$.  Each new n-dimensional block is a union of $N = (R_{j+1}/L)
... (R_l / L)$ blocks from the old decomposition.  More
generally, each interior p-face of the new decomposition is a
union of $N$ p-faces of the old decomposition. We let $B^+$ be
the complex generated by the interior faces of this coarser
decomposition.  Note that each interior face of $B^+$ has
dimension $p \ge l$ and dimensions $R_1 \times ... \times R_l
\times L \times ... \times L$.  (There are $p-l$ factors of $L$
in this formula.)

Any complex of cycles $C: B \rightarrow I_{rel}(S)$ can easily be
glued together to form a new complex of cycles $C^+: B^+
\rightarrow I_{rel}(S)$.  Suppose that $F$ is a p-face of $B^+$. 
As we observed above, $F$ is a union of p-faces from $B$: $F =
\sum_{i=1}^N F_i$, where $F_i$ is a face of $B$. Now we just
define $C^+(F) = \sum_{i=1}^N C(F_i)$.  The degree of $C$ and the
degree of $C^+$ are always the same.

In particular, $C_1^+$ is the glued-together version of $C_1$.
The volume of
$C_1^+(F^p)$ is at most $C(n) \delta N S_1 ... S_j S_j^{p-j}$. 
Plugging in the value of $N$, we see that the volume of
$C_1^+(F^p)$ is at most $C(n) \delta R_{j+1} ... R_l L^{-l+j} S_1
... S_j S_j^{p-j}$.  Finally, plugging in the value of $L$, we
see that the volume of $C_1^+(F^p)$ is at most

$$C(n,\delta) [\frac{R_1 ... R_j}{S_1 ... S_j}]^{\frac{l-k}{k-j}}
R_1 ... R_l S_j^{p-l}. \eqno{(V)}$$

Using the volume bound $(V)$ and the key lemma, we can now prove
estimates $(1)$ and $(2)$.  To prove inequality $(2)$, we set $l
= n$.  In this case $B^+$ consists of only one n-face, which is
the whole rectangle $R$.  According to the formula above,
$C_1^+(R)$ has volume at most $C(n) [\frac{R_1 ... R_j}{S_1 ...
S_j}]^{\frac{n-k}{k-j}} R_1 ... R_n$.  On the other hand, by the
Key Lemma, $C_1^+$ has degree $D$, and so $C_1^+(R)$ must have
volume at least $|D| S_1 ... S_n$. We conclude the following
inequality.

$$(R_1 ... R_j)^{\frac{n-k}{k-j}} R_1 ... R_n \ge c(n) |D| (S_1
... S_j)^{\frac{n-k}{k-j}} S_1 ... S_n.$$

This inequality is equivalent to $(2)$.

Next we prove inequality $(1)$ using Lemma 4.1.  Recall that
each interior face of $B^+$ has dimension $p \ge l$.  Since $C_1^+$
has degree $D \not= 0$, Lemma 4.1 guarantees that for some
dimension $p$, we can find an interior face $F^p$ so that
$C_1^+(F^p)$ has volume at least $c(n) S_1 ... S_p$.  On the other
hand, this same volume is bounded above by $(V)$.  Combining
these equations, we conclude the following.

$$C(n) [\frac{R_1 ... R_j}{S_1 ... S_j}]^{\frac{l-k}{k-j}}
R_1 ... R_l S_j^{p-l} \ge c(n) S_1 ... S_p.$$

Rearranging this inequality, we get the following.

$$[\frac{R_1 ... R_j}{S_1 ... S_j}]^{\frac{l-k}{k-j}}
R_1 ... R_l \ge c(n) S_1 ... S_p S_j^{-(p-l)} \ge c(n) S_1 ... S_l.$$

This proves inequality $(1)$.

We have now finished proving our main estimates except for the
proof of the key lemma which tells us that the degree of $C_1$ is
equal to $D$.

\vskip5pt

{\bf Gradually tightening chains}

\vskip5pt

\begin{key lemma} The degree of $C_1$ is equal to the degree of $\Phi$.
\end{key lemma}

To prove the lemma, we will need to construct some homotopies
between chain maps.  We use the following lemma, which
generalizes Lemma 4.1.

\begin{lemma} There is a constant $\epsilon(n) > 0$ so that the following
holds.  Suppose that $C_0$ and $C_1$ are two chain maps $X
\rightarrow I_{rel}(S)$.  Suppose that $C_0$ and $C_1$ agree on
the k-skeleton of $X$.  Suppose that for each p-face $F^p$ in $X$
of dimension $p \ge k+1$, the volumes $|C_0(F^p)|$ and
$|C_1(F^p)|$ are at most $\epsilon(n) S_1 ... S_p$.  Then $C_0$ and
$C_1$ are homotopic.
\end{lemma}

\proof We have to build a chain map $C: X \times [0,1] \rightarrow
I_{rel}(S)$, extending $C_0$ and $C_1$.  If $p \le k$, we define
$C(F^p \times [0,1])$ to be $0$.

We will prove inductively that we can extend $C$ to the
n-skeleton of $X \times [0,1]$ while preserving the
inequality $|C(F^p \times [0,1])| \le c(n) S_1 ... S_{p+1}$ 
for $p \le n-1$.

When we extend to the (p+1)-skeleton, we have to define $C(F^p
\times [0,1])$ for each p-face so that $\partial C(F^p \times [0,1]) =
C((\partial F^p) \times [0,1]) + C_1(F^p) - C_0(F^p)$.  By
induction, the right-hand side is a p-cycle in $S$ with volume at
most $c(n) S_1 ... S_p$.  According to our isoperimetric
inequality, we can fill this cycle with volume at most $c(n) S_1
... S_{p} S_{p} \le c(n) S_1 ... S_{p+1}$. 

Next we extend $C$ to the (n+1)-skeleton.  We have to define
$C(F^n \times [0,1])$.  We have already defined $C$ on $\partial(
F^n \times [0,1])$; it is an n-cycle with volume less than
$S_1 ... S_n$.  Hence it is an exact n-cycle, and we can choose
a filling for it.  We can then extend to the higher-dimensional
faces because $H_q(S, \partial S) = 0$ for all $q \ge n+1$.
\endproof

At first we might hope to apply this lemma to build a homotopy
from $C_0$ to $C_1$.  (Recall that $C_0$ and $C_1$ agree on the
k-skeleton of $B$.)  In general, this does not work, because the
volumes $|C_0(F^p)|$ may be too large.  Morally, the problem is
that in building $C_1$ we have suddenly tightened the chains into
a quite different position.  To build a homotopy, we want to
gradually tighten the chains so that at each step they move only
slightly.  Then we can use the lemma above to build a homotopy
between the small steps.

\vskip5pt

{\it Proof of key lemma}

\vskip5pt

Let $B_s$ be the division of $R$ into rectangular blocks with
dimensions $R_1 \times ... \times R_j \times 2^{-s} L \times ...
\times 2^{-s} L$.  The division $B_0$ is just $B$, and the
other $B_s$ are finer subdivisions of $B$.

Next we define chain maps $\Gamma_s: B_s \rightarrow I_{rel}(S)$ as
follows.  For each face $F^p$ in $B_s$ of dimension $p \le k$, we define 
$\Gamma_s(F)$ to be $\Phi(F \cap U)$.  Then we extend $\Gamma_s$ to
faces of dimension $p \ge k+1$ inductively, using the isoperimetric
inequality for rectangles at each step as in the construction of
$C_1$.  Because the constructions agree exactly, we may take
$\Gamma_0$ to be equal to $C_1$.

First we check that $\Gamma_{s+1}$ and $\Gamma_s$ have the same
degree.  We let $\Gamma_{s+1}^+: B_s \rightarrow I_{rel}(S)$ be
the glued version of $\Gamma_{s+1}$.  As in the previous gluing
construction, $\Gamma^+_{s+1}$ and $\Gamma_{s+1}$ have the same
degree.  We will use Lemma 6.1 to show that
$\Gamma_{s+1}^+$ and $\Gamma_s$ are homotopic.  By construction,
they have the same restriction to the k-skeleton of $B_s$.  By
the same argument that we used for $C_1$, $\Gamma_s(F^p)$ has
volume at most $\delta(n) S_1 ... S_j S_j^{p-j} \le \delta(n) S_1 ...
S_p$.  The same holds true for $\Gamma_{s+1}$ and hence for
$\Gamma^+_{s+1}$.  Applying Lemma 6.1, we see that
$\Gamma_s$ and $\Gamma_{s+1}$ have the same degree.

Let $\beta_s: B_s \rightarrow I_{rel}(S)$ be the chain map
sending a face $F^p$ to $\Phi(F^p \cap U)$ for every $p$.  The map
$\beta_s$ is analogous to $C_0$, and so it has degree $D$ for
every $s$.

If $p \ge k+1$, then the volume of $\beta_s(F^p)$ is at most
$|F^p|$, which is at most $R_1 ... R_j L^{p-j} 2^{-{p-j}s}$. 
Since $j < k$, we may choose $s$ sufficiently large so that for
each $p \ge k+1$, this volume is at most $c(n) S_1 ... S_p$.  
We now fix $s$ to be this sufficiently large value. The
two chain maps $\beta_s$ and $\Gamma_s$ agree on the k-skeleton
of $B_s$.  Because of our choice of $s$, the volume of
$\beta_s(F^p)$ is at most $c(n) S_1 ... S_p$ for each $p \ge
k+1$.  We checked above that the same inequality holds for
$\Gamma_s$.  According to Lemma 6.1, $\beta_s$ and
$\Gamma_s$ are homotopic and so have the same degree.

To summarize, the degree of $C_1$ is equal to the degree of
$\Gamma_0$, which is equal to the degree of $\Gamma_s$, which is
equal to the degree of $\beta_s$, which is equal to $D$.
\endproof

This concludes the proof of Estimates 1 and 2, and hence the
proof of Theorem 1.

\section{Constructing area-expanding embeddings}

In this section, we prove Theorem 2.

\begin{Theorem 2} There is a constant $C(n)$ so
that the following holds.

Suppose that the dimensions of $R$ and $S$ obey the following
inequalities for all $0 \le j < k \le l \le n$.

$$ (R_1 ... R_j)^{\frac{l-j}{k-j}} R_{j+1} ... R_l \ge C(n) (S_1
... S_j)^{\frac{l-j}{k-j}} S_{j+1} ... S_l. $$

Then there is a k-expanding embedding from $S$ into $R$.

\end{Theorem 2}

We will construct our embedding by composing a
k-expanding linear map and a simple folding map
analogous to the one in Figure 1.

If $R$ and $S$
are 2-dimensional rectangles with $R_1 > 3 S_1$ and $R_1 R_2 > 9 S_1
S_2$, then there is a 1-expanding embedding of $S$ into $R$.
This embedding is illustrated in Figure 1.

Next, let $a<b$ be integers between 1 and n.  If $R_i =
S_i$ except when $i$ is equal to $a$ or $b$ and
$R_a > 3 S_a$ and $R_a R_b > 9 S_a S_b$, then there is a
1-expanding embedding of $S$ into $R$.  This embedding is
the direct product of the folding map for the coordinates 
$a$ and $b$ and the identity in the other coordinates.

Composing these folding embeddings proves the following
lemma.

\begin{lemma} There is a constant $C(n)$ so that the following holds.
If $R_1 ... R_p > C(n) S_1 ... S_p$ for each $p$ between 1 and n, then
there is a 1-expanding embedding from $S$ into $R$.
\end{lemma}

The rest of the proof is just algebra, although it's rather tedious.
We put it in the form of a lemma.

\begin{lemma} Suppose that the dimensions of $R$ and $S$ obey the
following inequalities for all $0 \le j < k \le l \le n$.

$$ R_1 ... R_j (R_{j+1} ... R_l)^{\frac{k-j}{l-j}} \ge S_1
... S_j (S_{j+1} ... S_l)^{\frac{k-j}{l-j}}. \eqno{(In)}$$

Then there is a k-contracting linear diffeomorphism from $R$ to
a rectangle $T$ with $T_1 ... T_p \ge S_1 ... S_p$ for all $p$.
\end{lemma}

Given these lemmas, we finish the proof of Theorem 2. Under the
hypothesis of the theorem, Lemma 7.2 tells us that we can find a
k-contracting linear diffeomorphism from $R$ to $T$ where $T_1
... T_p \ge C(n) S_1 ... S_p$ for all $p$.  Then we use Lemma 7.1
to construct a 1-expanding embedding of $S$ into $T$.  Now we
turn to the proof of Lemma 7.2.

\proof If $S_1 ... S_p \le R_1 ... R_p$ for every p, then we take $T = R$
and we are done.  Let b be the smallest integer so that $S_1 ... S_b >
R_1 ... R_b$.  Because of all the inequalities in the hypothesis
of the lemma, we know that $b<k$.

We will define a sequence of linear diffeomorphisms $R = R(0)
\rightarrow R(1) \rightarrow ... \rightarrow R(c)$, for some
integer $c$ between 1 and $k-1$.  The diffeomorphism to $R(q)$ is
called $L_q$.  When $q$ is less than $c$, the rectangle $R(q)$ has
$R(q)_1 = ... =R(q)_{q+1}$.  The linear map $L_q$ increases each
$R(q-1)_i$ for $i$ between 1 and $q$ by a factor of $\lambda_q$ and
decreases every other $R(q-1)_i$ by a factor of
$\lambda_q^{-q/(k-q)}$, for some number $\lambda_q > 1$.  From
the last sentence, it follows that each $L_q$ is k-contracting. 
If $c$ is not bigger than $b$, then $R(c)_1 ... R(c)_b = S_1
... S_b$.  If $c$ is bigger than $b$, then $R(c)_1 ... R(c)_c = S_1
... S_c$.

Now we define the maps $L_q$.  It suffices to define $\lambda_q$. 
There is a maximum value of $\lambda_q$ which increases
$R(q-1)_q$ and decreases $R(q-1)_{q+1}$ until they meet.  If
there is a lesser value of $\lambda_q$ which makes $R(q)_1 ...
R(q)_m = S_1 ... S_m$, where $m$ is the maximum of $b$ and $q$, then
use that value and take $c=j$.  If not, use the maximal value. 
As we increase $q$, $R(q)_1 ... R(q)_b$ increases.  If $R(b)_1 ...
R(b)_b < S_1 ... S_b$, then $R(b)_1 ... R(b)_{b+1} < S_1 ...
S_{b+1}$, because $R(b)_1 = R(b)_{b+1}$.  More generally, for $q$
at least $b$, if $R(q)_1 ... R(q)_q < S_1 ... S_q$, then $R(q)_1
... R(q)_{q+1} < S_1 ... S_{q+1}$ also.

From the formula for the map $L_q$, it follows that $R(q)_1 ...
R(q)_k = R_1 ... R_k$ for every $q$, and by hypothesis $R_1 ... R_k
\ge S_1 ...S_k$.  Therefore, the above construction terminates
with $c$ less than or equal to $k-1$.

Recall that $m$ is the maximum of $b$ and $c$.  We have proven above that
$R(c)_1 ... R(c)_m = S_1 ... S_m$.  Moreover, for every $p$ less
than $m$, $R(c)_1 ... R(c)_p \ge S_1 ... S_p$.  If b is greater
than or equal to c, this follows because $R_1 ... R_p \ge S_1 ...
S_p$, and the definition of $L_j$ shows that $R(q)_1 ... R(q)_p
\ge R_1 ... R_p$ for every $p$ less than $k$.  If $c$ is greater than
$b$, this follows because $R(c)_1 = R(c)_m$ and $R(c)_1 ... R(c)_m
= S_1 ... S_m$.  In either case it is true.

The maps $L_q$ preserve many of the inequalities in $(In)$.
In particular, if $j \ge q$, then the following
equality holds.

$$R(q)_1 ... R(q)_j (R(q)_{j+1} ... R(q)_l)^{(k-j)/(l-j)} = R_1
... R_j (R_{j+1} ... R_l)^{(k-j)/(l-j)}.$$

Therefore, if $j \ge m$ then

$$R(c)_1 ... R(c)_j (R(c)_{j+1} ... R(c)_l)^{(k-j)/(l-j)} \ge S_1
... S_j (S_{j+1} ... S_l)^{(k-j)/(l-j)}.$$

Since $R(c)_1 ... R(c)_m = S_1 ... S_m$, we can divide the above
inequality on both sides, leaving the following inequality for
all $j \ge m$.

$$R(c)_{m+1} ... R(c)_j (R(c)_{j+1} ... R(c)_l)^{(k-j)/(l-j)} \ge
S_{m+1} ... S_j (S_{j+1} ... S_l)^{(k-j)/(l-j)}. \eqno{(*)}$$

At this point, we employ induction on the dimension of the rectangles.

We define $R'$ to be the (n-m)-directional rectangle with dimensions
$R(c)_{m+1} \times ... \times R_c(n)$, so that $R(c) = [0, R(c)_1] \times
... \times [0, R(c)_m] \times R'$.  We define $S' = S_{m+1} \times 
... \times S_n$.  We can rewrite $(*)$ in terms of $R'$
and $S'$.  To do this, let $k' = k-m$, $j' = j-m$ and $l' = l-m$.  Then
$(*)$ tells us that for any $j', l'$ in the ranges $0 \le j' < k' \le l'
\le n-m$, we have the following inequalities.

$$R'_1 ... R'_{j'} (R'_{j' + 1} ... R'_{l'})^{(k' - j')/ (l' -j')} \ge
S'_1 ... S'_{j'} (S'_{j' + 1} ... S'_{l'})^{(k'-j')/ (l' - j')}. \eqno{(*)'} $$

By induction on the dimension $n$, we can assume that there is a
$k'$-contracting linear diffeomorphism from $R'$ to some rectangle
$T'$ so that $T_1' ... T_p' \ge S_1' ... S_p'$ for any $1 \le p \le
n-m$.

We finally define T to be the rectangle with dimensions 
$R(c)_1 \times ... \times R(c)_m \times
T_1' \times ... \times T_{n-m}'$.  The direct product of the
(k-m)-contracting linear map from R' to T' with the identity map
is a k-contracting linear diffeomorphism from $R(c)$ to $T$.  Since
we already have a k-contracting linear map from $R$ to $R(c)$, we can
compose the two maps to get a k-contracting linear diffeomorphism 
from $R$ to
$T$.  Also, we already know that $T_1 ... T_p = R(c)_1 ... R(c)_p
\ge S_1 ... S_p$ when $p$ is less than or equal to $m$.  But for
larger $p$, $T_1 ... T_p = T_1 ... T_m T'_1 ... T'_{p-m} \ge S_1
... S_m S'_1 ... S'_{p-m} = S_1 ... S_p$.  \endproof

\section{Appendix: k-dilation and linear algebra}

In this section we record some basic facts about k-dilation
that follow from linear algebra.

If $L$ is a linear map from $\mathbb{R}^M$ to $\mathbb{R}^N$, 
then we can write $L$ in the form $0_1 D O_2$.  In this equation,
$O_2$ is an $M \times M$ orthogonal matrix, $O_1$ is an $N \times N$
orthogonal matrix, and $D$ is an $M \times N$ matrix which vanishes
off the diagonal and with all diagonal entries at least 0.
If we let $n$ 
be the minimum of $M$ and $N$, then the diagonal entries of $D$ are 
$0 \le s_1 \le ... \le s_n$.  The numbers $s_1, ...
, s_n$ are called the singular values of $L$.  The Lipschitz
constant of $L$ is the largest singular value $s_n$.
The k-dilation of $L$ is the product of the $k$ largest
singular values: $s_{n-k+1} ... s_n$.  Using this fact,
we prove some basic inequalities about k-dilation.

\begin{lemma} Suppose that $l > k$.  Then the following
inequality holds between the l-dilation and the k-dilation.

$$|\Lambda^l L|^{k/l} \le |\Lambda^k L|.$$

\end{lemma}

\proof Let $s_i$ denote the singular values of $L$.  Then the left-hand
side is $(s_{n-l+1} ... s_n)^{k/l}$.  This expression is less than
$(s_{n-k+1} ... s_n)$, which is the right-hand side. \endproof

\begin{corollary} Suppose that $\Phi$ is a map with k-dilation
$D(k)$ and l-dilation $D(l)$ with $l \ge k$.  Then
$D(l)^{k/l} \le D(k)$.
\end{corollary}

\proof Recall that $D(l)$ is the supremum of $| \Lambda^l d \Phi |$.  For
each point $x$, $| \Lambda^l d \Phi(x)|^{k/l} \le |\Lambda^k d \Phi(x)|$.
Passing to the supremum proves the corollary. \endproof

We include one other piece of linear algebra related to
k-dilation.  We don't use this result in our paper, but I think
it's worth knowing for context.  If $j < k$, a linear map with
k-dilation equal to 1 may have arbitrarily large j-dilation, but
it must pay for a large j-dilation by having a small l-dilation
for each $l > k$.  This tradeoff is described by the following
lemma.

\begin{lemma} Suppose that $j \le k \le l$ and that $L$ is a linear map.

$$\textrm{Then } |\Lambda^j L|^{\frac{l-k}{l-j}} |\Lambda^l
L|^{\frac{k-j}{l-j}} \le |\Lambda^k L|.$$

\end{lemma}

\proof The idea is to rewrite everything in terms of singular values.

$$|\Lambda^j L|^{l-k} |\Lambda^l L|^{k-j} = (s_{n-j+1} ...
s_n)^{l-k} (s_{n-l+1} ... s_n)^{k-j}$$

$$ = (s_{n-l+1} ... s_{n-j})^{k-j} (s_{n-j+1} ... s_n)^{l-j} \le
(s_{n-k+1} ... s_{n-j})^{l-j} (s_{n-j+1} ... s_n)^{l-j}$$

$$ = |\Lambda^k L|^{l-j}.$$

Taking $(l-j)^{th}$ roots of both sides finishes the proof. \endproof

I call this lemma the expansion/contraction inequality: a
k-contracting linear map may expand in some directions and has to
pay for it by contracting in others.  Unlike the last lemma, this
one has no direct analogue for non-linear maps.  A k-contracting
map $\Phi$ may have large j-dilation and l-dilation equal to 1. 
Suppose at one point $x$ that $|\Lambda^j d \Phi_x| = 10^6$.  It
follows that the j-dilation of $\Phi$ is at least $10^6$. It also
follows that the l-dilation of $\Phi$ at the point $x$ is small. 
But the l-dilation of $\Phi$ globally may still be 1 because at
some other point $y$, we may have $d\Phi_y$ equal to the
identity.

Nevertheless, the results in this paper can be viewed as an
analogue of the expansion/contraction inequality for nonlinear
maps.  For example, suppose that $\Phi$ is a degree 1
k-contracting map from $R$ to $S$.  Suppose that $R_1 ... R_j <<
S_1 ... S_j$.  Because of the sweepout lemma, the j-dilation of
$\Phi$ must be at least $c(n) S_1 ... S_j / R_1 ... R_j$.  If
$\Phi$ were linear, its l-dilation would then be bounded by a
small number coming from the expansion/contraction inequality. 
The actual l-dilation of $\Phi$ may be 1, but the map $\Phi$ can
be in some sense approximated by the complex of chains $C_1$
(defined in Section 6).  Up to a constant factor $C(n)$, the
volumes of the chains in $C_1$ obey the same bounds that would
follow if the l-dilation of $\Phi$ obeyed the
expansion/contraction inequality.

\section{Appendix 2: minor generalizations}

In this section, we discuss how far our results generalize to
shapes that are not rectangles.

First we briefly consider replacing $S$ by another shape.  We
note that all our arguments depended only on knowing the
isoperimetric profile of $S$.  Therefore, our methods should
adapt to give some estimates for any target where we can
estimate the isoperimetric profile.

Second we consider replacing $R$ by a more general shape.  Our
arguments apply to products of the form: $X^j \times Y^{l-j}
\times Z^{n-l}$, where $X$ and $Z$ may be any Riemannian
manifolds, but the middle factor $Y$ is still a rectangle.  In
this case, our estimate survives, reading as follows.

\begin{prop} Suppose that $U$ is an open set in $X \times Y
\times Z$, where $X^j$ and $Z^{n-l}$ are Riemannian manifolds
and $Y^{l-j}$ is a rectangle.  Suppose that $\Phi$ is a
k-contracting degree non-zero map from $U$ to an n-dimensional
rectangle $S$.  Then the volumes of $X$ and $Y$ are bounded below
by the following inequalities.

$$|X|^{\frac{l-j}{k-j}} |Y| \ge c(n) (S_1 ...
S_j)^{\frac{l-j}{k-j}} S_{j+1} ... S_l.$$

\end{prop}

If $X$ and $Z$ are oriented, $l=n$, and the degree of the map is
large, we also get an analogue of Estimate 2:
$|X|^{\frac{n-j}{k-j}} |Y| \ge c(n) |D| (S_1 ...
S_j)^{\frac{n-j}{k-j}} S_{j+1} ... S_n$.

\proof (sketch) Use the argument of the paper, cutting the domain
into pieces each a product of the form $X$ times a cube in $Y$
with side-length $L$ times a tiny simplex in $Z$. If the domain
is not orientable, use mod 2 chains instead of integral chains.
\endproof

The statement of the proposition would still make sense if we
allowed the middle factor $Y$ to be any manifold, but the
rectangular structure is used crucially in the proof, mostly when
we cut $Y$ into cubes.  I strongly believe that the estimate
above does not generalize to all Riemannian products $X \times Y
\times Z$.

The product structure can also be relaxed a little.  Suppose our
domain $A^n$ admits a map $\pi$ onto $Y^{l-j} \times Z^{n-l}$,
where as above $Y$ is a rectangle and $Z$ is any Riemannian
manifold.  Suppose that for any p-chain $C$ in $Y \times Z$, the
(p+j)-dimensional volume of $\pi^{-1}(C)$ is at most $V |C|$. 
Suppose that $U$ is an open set in $A$ admitting a k-contracting
map of non-zero degree to the n-dimensional rectangle $S$.  Then
our inequality again survives in the form $V^{\frac{l-j}{k-j}}
|Y| \ge c(n) (S_1 ... S_j)^{\frac{l-j}{k-j}} S_{j+1} ... S_l$. 
(And the analogue of Estimate 2 holds also.)

For example, we can replace $R$ by an ellipsoidal
metric on the n-sphere.  Define $E^n$ by the equation
$\sum_{i=0}^n (x_i / E_i)^2 = 1$.  Here $E$ is an ellipsoid
with principal axes $E_0 \le ... \le E_n$.  The manifold $E$ is
$C(n)$-bilipschitz to the double of the rectangle $[0, E_1]
\times ... \times [0, E_n]$.  So for any $j \ge 0$, there is a
map $\pi$ from $E$ to $[0, E_{j+1}] \times ... \times [0, E_n]$
which obeys the conditions of the last paragraph with $V \sim E_1
... E_j$.  Applying our generalized version of Estimates 1 and 2,
we get the following corollary.

\begin{cor} Suppose that $E$ and $E'$ are n-dimensional
ellipsoids with principal axes $E_0 \le ... \le E_n$ and $E'_0
\le ... \le E'_n$, and quotients $Q_i = E'_i / E_i$.  Suppose
that $\Phi$ is a map from $E$ to $E'$ with degree $D \not= 0$. 
Then the k-dilation of $\Phi$ is bounded below by the following
formulas.  First, if $0 \le j < k \le l \le n$,

$$dil_k(\Phi) \ge c(n) Q_1 ... Q_j (Q_{j+1} ...
Q_l)^{\frac{k-j}{l-j}}.$$

Second, if $0 \le j < k$, 

$$dil_k(\Phi) \ge c(n) |D|^{\frac{k-j}{n-j}} Q_1 ... Q_j (Q_{j+1}
... Q_n)^{\frac{k-j}{n-j}}.$$

\end{cor}

\end{document}